# ON THE MINIMAL TRAVEL TIME NEEDED TO COLLECT $n$ ITEMS ON A CIRCLE

By Nelly Litvak and Willem R. van Zwet

*University of Twente and University of Leiden*

Consider $n$ items located randomly on a circle of length 1. The locations of the items are assumed to be independent and uniformly distributed on $[0, 1)$. A picker starts at point 0 and has to collect all $n$ items by moving along the circle at unit speed in either direction. In this paper we study the minimal travel time of the picker. We obtain upper bounds and analyze the exact travel time distribution. Further, we derive closed-form limiting results when $n$ tends to infinity. We determine the behavior of the limiting distribution in a positive neighborhood of zero. The limiting random variable is closely related to exponential functionals associated with a Poisson process. These functionals occur in many areas and have been intensively studied in recent literature.

**1. Introduction.** This paper is devoted to the properties of the optimal route of the picker who has to collect $n$ items independently and uniformly distributed on a circle. By *optimal* route we mean the route providing the minimal travel time (see Figure 1). The problem has applications in performance analysis of carousel systems. A carousel is an automated storage and retrieval system which is widely used in modern warehouses. The system consists of a large number of shelves or drawers rotating in a closed loop in either direction. Orders are represented by a list of items. The list specifies the type and retrieval quantity of each item. The picker has a fixed position in front of the carousel, which rotates the required items to the picker. In this paper we study the minimal travel (rotation) time of the carousel while picking one order of $n$ items, the locations of which are assumed to be independent and uniformly distributed on the carousel.

Let $U_0 = 0$ be the picker's starting point and, for $i = 1, 2, \ldots, n$, let the random variable $U_i$ denote the position of the $i$th item. The random vari-

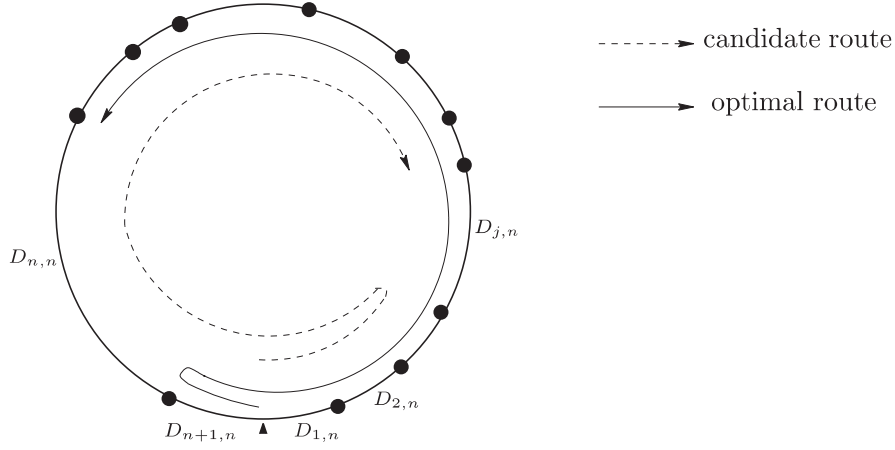

Fig. 1. *Minimal travel time on a circle.*

ables $U_1, U_2, \ldots, U_n$ are independent and uniformly distributed on $[0,1)$. Set $U_{n+1} = 1$. Let

$$0 = U_{0:n} < U_{1:n} < \cdots < U_{n:n} < U_{n+1:n} = 1$$

denote the ordered $U_0, U_1, \ldots, U_{n+1}$. Then the picker's starting point and the positions of the $n$ items partition the circle into $n+1$ uniform spacings

$$D_{i,n} = U_{i:n} - U_{i-1:n}, \qquad 1 \leq i \leq n+1.$$

Let $X_1, X_2, \ldots$ be independent exponential random variables with mean 1 and write

$$S_0 = 0, \qquad S_i = \sum_{j=1}^{i} X_j, \qquad i \geq 1.$$

It is well known that [cf. Pyke (1965)]

(1.1) $(D_{1,n}, D_{2,n}, \ldots, D_{n+1,n}) \stackrel{d}{=} (X_1/S_{n+1}, X_2/S_{n+1}, \ldots, X_{n+1}/S_{n+1}),$

that is, the spacings are distributed as normalized exponentials. According to Pyke (1965), this construction is useful "to show that an ordering of uniform spacings may be considered as an ordering of the exponential random variables."

Now, let $T_n$ be the minimal travel time. We explore $T_n$ in terms of the uniform $(n+1)$-spacings $D_{1,n}, D_{2,n}, \ldots, D_{n+1,n}$. For $n = 1$, the problem is trivial. The picker just chooses the shorter distance from the starting point to the item, and thus the travel time $T_1$ is distributed as $(1/2)D_{1,1}$ (a normalized minimum of two exponentials). For $n = 2$, one can easily verify that the optimal route is guaranteed by the nearest item heuristic where the next item to be picked is always the nearest one. The travel time distribution for this heuristic was obtained by Litvak and Adan (2001). It follows from their result that $T_2$ is distributed as $(1/2)D_{1,2} + (3/4)D_{2,2}$. For $n \geq 3$, the problem becomes much more difficult.

A crucial and simple observation made by many authors [see, e.g., Bartholdi and Platzman (1986)] is that the optimal route admits at most one turn. Obviously, it is never optimal to cover the same segment of the circle more than twice. Thus, in general, $T_n$ can be expressed as

$$
(1.2) \quad T_n = 1 - \max\Big\{\max_{1 \leq j \leq n}\{D_{j,n} - U_{j-1:n}\}, \max_{1 \leq j \leq n}\{D_{n+2-j,n} - (1 - U_{n+2-j:n})\}\Big\}.
$$

This formula is easy to understand by means of Figure 1. Clearly, for $j = 1, 2, \ldots, n$, the term $D_{j,n} - U_{j-1:n}$ is the gain in travel time (compared to one full rotation) obtained by skipping the spacing $D_{j,n}$ and going back instead (ending in a clockwise direction). The same can be said about $D_{n+2-j,n} - (1 - U_{n+2-j:n})$, but here the picker ends in a counterclockwise direction. Under the optimal strategy, the picker chooses the largest possible gain.

Let $T_n^{(m)}$ be the travel time under so-called $m$-step strategies: the picker chooses the shortest route among $2(m+1)$ candidate routes that change direction at most once (as does the optimal route) and only do so after collecting no more than $m$ items. It was proved by Litvak and Adan (2002) that for $2m < n$,

$$
(1.3) \quad \begin{aligned}
T_n^{(m)} &= 1 - \max\Big\{\max_{1 \leq j \leq m+1}\{D_{j,n} - U_{j-1:n}\}, \\
&\qquad\qquad\qquad \max_{1 \leq j \leq m+1}\{D_{n+2-j,n} - (1 - U_{n+2-j:n})\}\Big\} \\
&\stackrel{d}{=} 1 - \frac{1}{S_{n+1}}\max\Big\{\max_{1 \leq j \leq m+1}\{X_j - S_{j-1}\}, \\
&\qquad\qquad\qquad \max_{1 \leq j \leq m+1}\{X_{n+2-j} - (S_{n+1} - S_{n+2-j})\}\Big\} \\
&\stackrel{d}{=} 1 - \max\Big\{\sum_{j=1}^{m+1}\frac{1}{2^j - 1}D_{j,n}, \sum_{j=1}^{m+1}\frac{1}{2^j - 1}D_{n+2-j,n}\Big\}.
\end{aligned}
$$

Formula (1.3) follows from the following curious property of exponential random variables obtained by Litvak (2001).

LEMMA 1.1. *For any* $m = 0, 1, \ldots$ *and* $0 < q < 1$,

$$
(1.4) \quad \max_{1 \leq j \leq m+1}\{X_j - (q^{-1} - 1)S_{j-1}\} \stackrel{d}{=} (q^{-1} - 1)\sum_{j=1}^{m+1} q^j(1 - q^j)^{-1}X_j.
$$

The proof also implies that for any $m = 0, 1, \ldots, n$,

$$(1.5) \quad \max_{1 \leq j \leq m+1} \{D_{j,n} - (q^{-1} - 1)U_{j-1:n}\} \stackrel{d}{=} (q^{-1} - 1) \sum_{j=1}^{m+1} q^j (1 - q^j)^{-1} D_{j,n}.$$

If $2m < n$, then the two internal maxima in the third expression of (1.3) are independent and (1.4) can be used to rewrite each of them separately. Moreover, the same argument applies for the normalized exponentials yielding (1.3). If $2m \geq n$, then the two internal maxima become dependent and the argument fails.

In fact, the optimal strategy is the $m$-step strategy with $m = n - 1$. Intuitively, it is clear, however, that with high probability the optimal route has only a few steps before a turn. That is, the $m$-step strategy often prescribes the optimal picking sequence even when $m$ is relatively small. It was shown by Litvak and Adan (2002) that already for $m = 2$, the $m$-step strategy is quite close to optimal and, on average, outperforms the nearest item heuristic.

Let $K_n^{(m)}$ and $K_n$ denote a number of items collected before a turn under the $m$-step strategy and the optimal strategy, respectively. If there is no turn, these numbers are set equal to zero. It was proved by Litvak and Adan (2002) that:

(i) $T_n^{(m)}$ and $K_n^{(m)}$ are independent random variables; (ii) for any $k = 0, 1, \ldots, m$ and $2m < n$,

$$(1.6) \quad \begin{aligned} \mathbb{P}(K_n^{(m)} = k) &= \mathbb{P}\bigg(\bigg[\arg\max_{1 \leq j \leq m+1}\{D_{j,n} - U_{j-1:n}\} = k+1\bigg]\bigg) \\ &= \mathbb{P}\bigg(\bigg[\arg\max_{1 \leq j \leq m+1}\{X_j - S_{j-1}\} = k+1\bigg]\bigg) \\ &= \frac{1}{2^{k+1} - 2^{k-m}}; \end{aligned}$$

(iii) for any $k = 0, 1, \ldots, n - 2$,

$$(1.7) \quad \mathbb{P}(K_n > k) < 1/2^k.$$

The last estimate is helpful in the analysis of the limiting properties of the optimal route. For example, it was proved by Litvak and Adan (2002) that for any fixed $k = 0, 1, \ldots,$

$$(1.8) \quad \lim_{n \to \infty} \mathbb{P}(K_n = k) = 1/2^{k+1}.$$

Indeed, observe that for $k = 0, 1, \ldots, m$,

$$\mathbb{P}(K_n^{(m)} = k) - \mathbb{P}(K_n > m) \leq \mathbb{P}(K_n = k) \leq \mathbb{P}(K_n^{(m)} = k).$$

Now, let $m$ and $n$ go to infinity in such a way that the inequality $2m < n$ is always satisfied. Then (1.8) follows readily from (1.6) and (1.7).

In this paper we first derive simple upper bounds for the minimal travel time. Then we analyse the distribution of $T_n$. Further, we obtain the limiting behavior of $T_n$ when $n$ tends to infinity.

**2. Upper bounds.** Let $T_n$ be the minimal travel time needed to collect $n$ items independently and uniformly distributed on a circle of length 1. The following lemma gives an upper bound that holds for *any* realization of the random items' locations.

LEMMA 2.1. *For any $n \geq 1$, the travel time $T_n$ never exceeds $1 - \alpha_{n+1}$, where*

$$\alpha_{n+1} = \frac{1}{2^{m+1} + 2^m - 2}, \quad n = 2m;$$

$$\alpha_{n+1} = \frac{1}{2 \cdot 2^{m+1} - 2}, \quad n = 2m+1.$$

*This upper bound is tight.*

PROOF. Assume that $n = 2m+1$. For $n = 2m$ the proof is similar. The positions of the items plus the picker's starting point partition the circle into $n+1$ spacings with lengths $d_1, d_2, \ldots, d_{n+1}$. Note that for any collection $d_1, d_2, \ldots, d_{n+1} \geq 0$ there exists a number $j = 1, 2, \ldots, m+1$ such that either (i) $d_j \geq 2^{j-1}\alpha_{n+1}$, $d_l < 2^{l-1}\alpha_{n+1}$, $l = 1, 2, \ldots, j-1$, or (ii) $d_{n+2-j} \geq 2^{j-1}\alpha_{n+1}$, $d_{n+2-l} < 2^{l-1}\alpha_{n+1}$, $l = 1, 2, \ldots, j-1$. This follows since

$$2 \sum_{j=1}^{m+1} 2^{j-1}\alpha_{n+1} = d_1 + d_2 + \cdots + d_{n+1} = 1.$$

Without loss of generality assume (i). Then the route that skips the spacing $d_j$ and goes back instead has length

$$1 - d_j + d_1 + d_2 + \cdots + d_{j-1} \leq 1 - \alpha_{n+1},$$

and its length must be greater or equal than $T_n$. This proves the upper bound.

To show the tightness we just put $d_j = d_{n+2-j} = 2^{j-1}\alpha_{n+1}$, $j = 1, 2, \ldots, m+1$. In this case the travel time under the optimal strategy equals $1 - \alpha_{n+1}$. □

Let us now consider the following approximation of $T_n$ in (1.2),

$$T_n^0 \stackrel{d}{=} 1 - \frac{1}{S_{n+1}} \max \Big\{ \max_{1 \leq j \leq m+1} \{X_j - S_{j-1}\},$$

$$\max_{1 \leq j \leq m'+1} \{X_{n+2-j} - (S_{n+1} - S_{n+2-j})\} \Big\},$$

where $m = m' = (n-1)/2$ if $n$ is odd and $m = m' + 1 = n/2$ if $n$ is even. In both cases we have $m + m' = n - 1$ so that the $X_j$'s from the first internal maximum are not involved in the second internal maximum. That is, the two internal maxima are independent, and we can apply Lemma 1.1 to each of these separately to arrive at

$$(2.1) \qquad T_n^0 \stackrel{d}{=} 1 - \max\left\{\sum_{j=1}^{m+1} \frac{1}{2^j-1}D_{j,n}, \sum_{j=1}^{m'+1}\frac{1}{2^j-1}D_{n+2-j,n}\right\}.$$

Clearly, $T_n^0$ gives a tight stochastic upper bound for $T_n$. In fact, $T_n^0$ and $T_n$ differ with probability of order $2^{-n/2}$ according to (1.7). It was shown by Litvak and Adan (2002) that $T_n^0$ is stochastically larger than the weighted sum

$$T_n^* = \sum_{j=2}^{n+1}(1-\alpha_j)D_{j,n}.$$

Straightforward estimation of the expected difference between $T_n^*$ and $T_n^0$ yields

$$\mathbb{E}(T_n^0 - T_n^*) < 0.09\mathbb{E}(D_{1,n}) = \frac{0.09}{n+1}.$$

Thus,

$$(2.2) \qquad \mathbb{E}(T_n) < \mathbb{E}(T_n^0) < \frac{1}{n+1}\sum_{j=2}^{n+1}(1-\alpha_j) + \frac{0.09}{n+1}.$$

In Table 1 (see Section 4), we compare the mean travel time $\mathbb{E}(T_n)$ obtained by simulation with upper estimate (2.2) and approximation (4.8), which follows from the limiting results in Section 4. The results prove that the bound (2.2) is quite sharp. For larger $n$, however, (4.8) gives a slightly better approximation.

**3. The minimal travel time distribution.** In this section we produce an explicit expression for $\mathbb{P}(T_n \geq 1-t)$. First, note that it is never optimal to turn after covering half of a circle. Now, consider the events

$$A_{n,k}(u,v) = [U_{k:n} = u < 1/2 < 1-v = U_{k+1:n}],$$
$$0 \leq u, v < 1/2, k = 0, 1, \ldots, n.$$

For $k = 2, 3, \ldots, n-2$, the joint distribution of $U_{1:n}, \ldots, U_{k-1:n}, 1-U_{k+2:n}, \ldots, 1-U_{n:n}$ given $A_{n,k}(u,v)$ is that of

$$uU_{1:k-1}, \ldots, uU_{k-1:k-1}, \qquad vV_{n-k-1:n-k-1}, \ldots, vV_{1:n-k-1},$$

where **U** and **V** are independent vectors of uniform order statistics. As the event $[T_n \geq 1-t]$ implies $1-v-u-u\wedge v \leq t$, we have for $k = 2, 3, \ldots, n-2$,

$$
\begin{aligned}
&\mathbb{P}(T_n \geq 1-t|A_{n,k}(u,v)) \\
&= \mathbb{P}\Big(\max_{1\leq j\leq k}\{(U_{j:k-1} - U_{j-1:k-1}) - U_{j-1:k-1}\} \leq t/u, \\
&\qquad \max_{1\leq j\leq n-k}\{(V_{j:n-k-1} - V_{j-1:n-k-1}) - V_{j-1:n-k-1}\} \leq t/v\Big) \\
&\quad \times \mathbf{1}_{[1-v-u-u\wedge v\leq t]} \\
&= P_{k-1}(t/u)P_{n-k-1}(t/v)\mathbf{1}_{[1-v-u-u\wedge v\leq t]}.
\end{aligned}
\tag{3.1}
$$

Here $u \wedge v = \min\{u, v\}$ denotes the smaller of $u$ and $v$ and

$$
P_m(t) = \mathbb{P}\Big(\max_{1\leq j\leq m+1}\{D_{j,m} - U_{j-1:m}\} \leq t\Big), \qquad m = 1, 2, \ldots; t \geq 0.
\tag{3.2}
$$

One readily verifies that the final expression in (3.1) continues to hold for $k = 1$ and $k = n-1$, provided we define

$$
P_0(t) = \mathbf{1}_{[t>1]}, \qquad t \geq 0.
\tag{3.3}
$$

For $k = 0$ and $k = n$, we find

$$
\mathbb{P}(T_n \geq 1-t|A_{n,0}(0,v)) = P_{n-1}(t/v)\mathbf{1}_{[1-v\leq t]},
$$
$$
\mathbb{P}(T_n \geq 1-t|A_{n,n}(u,0)) = P_{n-1}(t/u)\mathbf{1}_{[1-u\leq t]}.
$$

It follows that

$$
\begin{aligned}
&\mathbb{P}(T_n \geq 1-t) \\
&= \int_0^{1/2}\int_0^{1/2}\sum_{k=1}^{n-1}\binom{n}{k}ku^{k-1}(n-k)v^{n-k-1} \\
&\qquad \times P_{k-1}(t/u)P_{n-k-1}(t/v)\mathbf{1}_{[1-v-u-u\wedge v\leq t]}\,du\,dv \\
&\quad + 2\cdot\mathbf{1}_{[t>1/2]}\int_{u=1-t}^{1/2}nu^{n-1}P_{n-1}(t/u)\,du.
\end{aligned}
\tag{3.4}
$$

Formula (1.5) and Theorem 2 of Ali and Obaidullah (1982) imply an expression for $P_m(t)$. Writing

$$c_j = (2^j - 1)^{-1}, \qquad j = 1, 2, \ldots,$$

and $x_+ = \max\{x, 0\}$ for the positive part of a number $x$, we find that for $m = 1, 2, \ldots$,

$$
P_m(t) = \mathbb{P}\bigg(\sum_{j=1}^{m+1}c_jD_{j,m} \leq t\bigg) = \sum_{j=1}^{m+1}\{(t-c_j)_+\}^m\prod_{\substack{l=1,\\l\neq j}}^{m+1}(c_l - c_j)^{-1}.
\tag{3.5}
$$

The last expression is also valid for $m=0$. Of course, for $t>1$, the terms in (3.5) sum to 1.

Alternatively, one can determine $P_m(t)$, recursively. Conditioning on $U_{1:m}$, we find for $m=2,3,\ldots,$

$$P\left(\max_{1\leq i\leq m+1}\{D_{i,m}-U_{i-1:m}\}<t\,\Big|\,U_{1:m}=u\right)$$
$$=P\left((1-u)\max_{1\leq i\leq m}\{D_{i,m-1}-U_{i-1:m-1}\}-u<t\right)\mathbf{1}_{[u\leq t]}$$
$$=P_{m-1}\left(\frac{t+u}{1-u}\right)\mathbf{1}_{[u\leq t]}.$$

This yields the recursive equation

$$(3.6)\qquad P_m(t)=\int_0^t m(1-u)^{m-1}P_{m-1}\left(\frac{t+u}{1-u}\right)du,$$

which is valid for $m=1,2,\ldots.$

We can now find the distribution of $T_n$ by substituting (3.3) and either (3.5) or (3.6) in (3.4) and integrating. One obtains, for example, for $t\geq 0$,

$$P_1(t)=\tfrac{1}{2}(3t-1)_+-\tfrac{3}{2}(t-1)_+,$$
$$P_2(t)=\tfrac{1}{8}\{(7t-1)_+\}^2-\tfrac{7}{8}\{(3t-1)_+\}^2+\tfrac{7}{4}\{(t-1)_+\}^2$$

and for $0\leq t\leq 1$,

$$\mathbb{P}(T_1\geq 1-t)=(2t-1)_+,$$
$$\mathbb{P}(T_2\geq 1-t)=\tfrac{1}{3}\{(4t-1)_+\}^2-2\{(2t-1)_+\}^2,$$
$$\mathbb{P}(T_3\geq 1-t)=\tfrac{1}{4}\{(6t-1)_+\}^3-\tfrac{41}{36}\{(4t-1)_+\}^3$$
$$\qquad\qquad-\tfrac{1}{4}\{(4t-1)_+\}^2+\tfrac{11}{4}\{(2t-1)_+\}^3.$$

Although the general structure of these functions is fairly easy to understand, it seems quite useless to provide explicit expressions for $\mathbb{P}(T_n\geq 1-t)$ for much larger values of $n$. Instead, we study their asymptotic behavior in Section 4.

**4. Limiting results.** In this section we shall obtain the limiting distribution of $(n+1)(1-T_n)$. First of all, let us consider the limiting behavior of

$$P_m(t/(m+1))=\mathbb{P}\left((m+1)\sum_{j=1}^{m+1}\frac{1}{2^j-1}D_{j,m}<t\right).$$

THEOREM 4.1. *Let $X_1, X_2, \ldots$ be independent exponential random variables with mean 1. Then*

$$(4.1) \qquad (m+1) \sum_{j=1}^{m+1} \frac{1}{2^j - 1} D_{j,m} \xrightarrow{d} \sum_{j=1}^{\infty} \frac{1}{2^j - 1} X_j,$$

*and the limiting distribution is given by*

$$(4.2) \qquad \begin{aligned} P(t) &= \lim_{m \to \infty} P_m(t/(m+1)) \\ &= 1 - \sum_{j=1}^{\infty} (-1)^{j-1} 2^j \exp\{-(2^j - 1)t\} \prod_{l=1}^{j} \frac{1}{2^l - 1}. \end{aligned}$$

*The distribution function $P$ satisfies the integral equation*

$$(4.3) \qquad e^{-t} P(t) = \int_t^{2t} e^{-u} P(u) \, du.$$

PROOF. The argument essentially repeats the proof of Theorem 4 of Litvak and Adan (2001). Define

$$(4.4) \qquad J_m = \sum_{j=1}^{m+1} \frac{1}{2^j - 1} X_j, \qquad J = \sum_{j=1}^{\infty} \frac{1}{2^j - 1} X_j.$$

By the monotone convergence theorem, $\mathbb{E}(J) = \lim_{m \to \infty} \mathbb{E}(J_m) < \infty$. In particular, it implies $\mathbb{P}(J < \infty) = 1$.

Now, using (1.1), we write

$$(m+1) \sum_{j=1}^{m+1} \frac{1}{2^j - 1} D_{j,m} \stackrel{d}{=} \frac{(m+1) J_m}{S_{m+1}}.$$

By definition, the sequence $\{J_m\}$ converges a.s. to $J$. The strong law of large numbers implies that the sequence $\{(m+1)/S_{m+1}\}$ converges a.s. to 1. Thus, $\{(m+1)J_m/S_{m+1}\}$ converges a.s. to $J$ which immediately gives (4.1).

The distribution $P$ of $J$ can be obtained via inversion of its Laplace–Stieltjes transform

$$\varphi(s) = \mathbb{E}(\exp(-sJ)) = \prod_{j=1}^{\infty} \frac{2^j - 1}{2^j - 1 + s}.$$

One can expand $\varphi(s)$ in rational fractions of $s$ and obtain

$$(4.5) \qquad \varphi(s) = \sum_{j=1}^{\infty} \frac{(-1)^{j-1} 2^j}{2^j - 1 + s} \prod_{l=1}^{j-1} \frac{1}{2^l - 1}.$$

Here, in order to write the formula for the residues of $\varphi(s)$, one can apply well-known expressions from so-called $q$-calculus [see, e.g., Gasper and Rahman (1990)], but in our case it is not difficult to verify this formula directly. Inversion of (4.5) yields (4.2).

Finally, we use (3.6) and the dominated convergence theorem to obtain

$$
\begin{aligned}
P(t) &= \lim_{m\to\infty} P_m(t/(m+1)) \\
&= \lim_{m\to\infty} \int_0^{t/(m+1)} m(1-u)^{m-1} P_{m-1}\left(\frac{t/(m+1)+u}{1-u}\right) du \\
&= \int_0^t e^{-u} P(t+u)\, du \\
&= e^t \int_t^{2t} e^{-u} P(u)\, du,
\end{aligned}
$$

which proves (4.3). $\square$

Obviously, we also have convergence of moments. For the $k$th moment of $P$, (4.2) yields

$$\mathbb{E}(J^k) = \int_0^\infty t^k\, dP(t) = k! \sum_{j=1}^\infty (-1)^{j-1} \frac{2^j}{(2^j-1)^k} \prod_{l=1}^j \frac{1}{2^l-1}.$$

Alternatively, one can directly use (4.4) to find a simple expression for cumulants $\kappa_\nu$, $\nu \geq 1$, of $P$. It is immediate that

$$\mathbb{E}(J) = \kappa_1 = \sum_{j=1}^\infty (2^j-1)^{-1},$$

$$\mathrm{Var}(J) = \kappa_2 = \sum_{j=1}^\infty (2^j-1)^{-2}.$$

Furthermore,

$$\log(\mathbb{E}\exp(itJ)) = -\sum_{j=1}^\infty \log(1-(2^j-1)^{-1} it) = \sum_{j=1}^\infty \sum_{\nu=1}^\infty \frac{(it)^\nu}{\nu(2^j-1)^\nu},$$

where $i$ is the imaginary unit. Since $\log(\mathbb{E}\exp(itJ)) = \sum_{\nu=1}^\infty \kappa_\nu (it)^\nu (\nu!)^{-1}$, it follows that

$$\kappa_\nu = (\nu-1)! \sum_{j=1}^\infty (2^j-1)^{-\nu}, \qquad \nu \geq 1.$$

The distribution function $P$ on $[0,\infty)$ has the remarkable property that it is infinitely often differentiable and that all of its derivatives $P^{(k)}$ vanish at the origin. This is most easily seen by differentiating (4.3), but one may also use (4.2) to show analytically that $P^{(k)}(0) = 0$ for all $k = 1, 2, \ldots$. It follows that $P$ is not analytic at the origin. The series (4.2) diverges for all $t < 0$ and, hence, $P$ cannot be represented by its Taylor series around $t = 0$.

Now repeating the argument from the proof of Theorem 4.1, one can show that

$$(n+1)(1 - T_n^0) \xrightarrow{d} \max\left\{\sum_{j=1}^{\infty} \frac{1}{2^j - 1} X_j, \sum_{j=1}^{\infty} \frac{1}{2^j - 1} X_j'\right\},$$

where $X_1, X_2, \ldots, X_1', X_2', \ldots$ are independent exponentials with mean 1. Since the two sums in the maximum are independent and (1.7) ensures that

$$\mathbb{P}(T_n \neq T_n^0) < 2^{-(n-2)/2},$$

we have proved the following statement.

THEOREM 4.2. *Let $X_1, X_2, \ldots, X_1', X_2', \ldots$ be independent exponential random variables with mean 1. Then*

(4.6) $$(n+1)(1 - T_n) \xrightarrow{d} \max\left\{\sum_{j=1}^{\infty} \frac{1}{2^j - 1} X_j, \sum_{j=1}^{\infty} \frac{1}{2^j - 1} X_j'\right\},$$

*and the limiting distribution is*

(4.7) $$\lim_{n \to \infty} \mathbb{P}(T_n > 1 - t/(n+1)) = [P(t)]^2,$$

*where $P(t)$ is given by (4.2).*

Again we have moment convergence and for the $k$th moment we find

$$\lim_{n \to \infty} \mathbb{E}[(n+1)(1 - T_n)]^k$$

$$= 2k! \sum_{j=1}^{\infty} (-1)^{j-1} \frac{2^j}{(2^j - 1)^k} \prod_{l=1}^{j} \frac{1}{2^l - 1}$$

$$- 2k! \sum_{j=1}^{\infty} \sum_{i=1}^{\infty} (-1)^{i+j} \frac{2^{i+j}}{(2^i + 2^j - 2)^{k+1}} \prod_{l=1}^{j} \frac{1}{2^l - 1} \prod_{r=1}^{i-1} \frac{1}{2^r - 1}.$$

An equivalent expression for the expectation can be obtained as

$$\lim_{n \to \infty} \mathbb{E}[(n+1)(1 - T_n)]$$

$$= \int_0^{\infty} (1 - [P(t)]^2) \, dt$$

$$= 2\sum_{j=1}^{\infty} \frac{1}{2^j - 1} - \sum_{i=1}^{\infty}\sum_{j=1}^{\infty} (-1)^{i+j} \frac{2^{i+j}}{2^i + 2^j - 2} \prod_{l=1}^{j} \frac{1}{2^l - 1} \prod_{r=1}^{i} \frac{1}{2^r - 1}$$
$$\approx 2.1578.$$

For large $n$ we therefore have the estimate

(4.8) $$\mathbb{E}(T_n) \approx 1 - \frac{2.1578}{n+1}.$$

In Table 1 we compare the mean travel time obtained by simulation with upper estimate (2.2) (see Section 2) and approximation (4.8). We see that both approximations are quite sharp, but (4.8) performs somewhat better. It is no surprise that both (2.2) and (4.8) are close to $\mathbb{E}(T_n)$ for large $n$ since all three quantities converge to 1 as $n \to \infty$. What is encouraging is that, already for $n = 30$, both approximations of $(n+1)(1 - \mathbb{E}(T_n))$ are very good. That (4.8) yields a better approximation of $\mathbb{E}(T_n)$ than (2.2) is to be expected since it is asymptotically correct up to and including order $n^{-1}$, whereas (2.2) has a slight asymptotic error of about $+0.006/(n+1)$. After all, (2.2) was derived as an upper bound.

Table 1
*Estimation of the mean travel time under the optimal strategy*

| $n$ | 3 | 5 | 10 | 15 | 20 | 30 |
|---|---|---|---|---|---|---|
| $\mathbb{E}(T_n)$ | 0.5262 | 0.6591 | 0.8052 | 0.8653 | 0.8972 | 0.9304 |
| $\mathbb{E}[(n+1)(1-T_n)]$ | 1.8952 | 2.0454 | 2.1423 | 2.1548 | 2.1592 | 2.1572 |
| Upper estimate (2.2) | 0.5433 | 0.6670 | 0.8068 | 0.8658 | 0.8976 | 0.9306 |
| $(n+1)$[1-upper estimate(2.2)] | 1.8268 | 1.9980 | 2.1252 | 2.1472 | 2.1504 | 2.1514 |
| Approximation (4.8) | 0.4605 | 0.6404 | 0.8038 | 0.8651 | 0.8972 | 0.9304 |
| $(n+1)$[1-approximation (4.8)] | 2.1578 | 2.1578 | 2.1578 | 2.1578 | 2.1578 | 2.1578 |

**5. Asymptotic behavior in the neighborhood of zero.** In this section we study the behavior of $P(t)$ as $t \to +0$. So far we have found only that $P$ has vanishing derivatives at the origin and can not be expanded in a Taylor expansion around $t = 0$. We shall, therefore, have to attack this problem in a different manner.

Let $X_1, X_2, \ldots$ be independent exponential random variables with mean 1, let

$$c_j = (2^j - 1)^{-1}, \qquad j = 1, 2, \ldots,$$

and define

$$J = \sum_{j=1}^{\infty} c_j X_j.$$

We want to determine the behavior of

$$P(t) = \mathbb{P}(J \leq t)$$

for small positive values of $t$. In principle this problem is solved in Theorem 3.2 of Davis and Resnick (1991), but we need to do a substantial amount of analysis to make their result explicit, even in our relatively simple case.

In our case, the distribution function $F(x) = P(X_1 < x) = 1 - \exp\{-x\}$ and the density $f(x) = \exp\{-x\}$ are regularly varying at 0 with index $\alpha = 1$ and $\alpha - 1 = 0$, respectively. The $c_j$'s are positive and nonincreasing, their sum converges and for every $\theta \in (0,1)$,

$$\theta^n \sum_{j=1}^{\infty} \{c_j^2/c_n^2\} \mathbf{1}_{[j \geq \theta^{-n}]} = \theta^n \sum_{j=1}^{\infty} \{(2^n-1)/(2^j-1)\}^2 \mathbf{1}_{[j \geq \theta^{-n}]} \to 0$$

as $n \to \infty$. The density $f$ satisfies

$$\int_0^{\infty} e^{-2\lambda x} f^2(x) \, dx = 1/\{2(1+\lambda)\} \qquad \text{for } \lambda > 0.$$

Hence, we have verified the assumptions of Theorem 3.2 of Davis and Resnick (1991) in our case. The theorem states that

(5.1) $\qquad P(m_\lambda) \sim \exp\{\lambda m_\lambda\} \varphi_J(\lambda)/(\lambda S_\lambda \sqrt{2\pi}) \qquad \text{as } \lambda \to \infty.$

Here

$$m_\lambda = \sum_{j=1}^{\infty} \frac{c_j}{1+\lambda c_j} = \sum_{j=1}^{\infty} \frac{1}{2^j - 1 + \lambda},$$

$$\varphi_J(\lambda) = \prod_{j=1}^{\infty} \frac{1}{1+\lambda c_j} = \prod_{j=1}^{\infty} \frac{2^j - 1}{2^j - 1 + \lambda}$$

and

$$S_\lambda^2 = \sum_{j=1}^{\infty} \frac{c_j^2}{(1+\lambda c_j)^2} = \sum_{j=1}^{\infty} \frac{1}{(2^j - 1 + \lambda)^2}.$$

We obviously have to study the behavior of these quantities as $\lambda \to \infty$ and, hence, $m_\lambda \to 0$. It is easier to deal with integrals than sums. For $k = 1, 2, \ldots$ and $\lambda \to \infty$, we have

$$0 \leq \int_0^{\infty} (2^x - 1 + \lambda)^{-k} \, dx - \sum_{j=1}^{\infty} (2^j - 1 + \lambda)^{-k}$$

$$\leq \sum_{j=0}^{\infty}(2^j - 1 + \lambda)^{-k} - \sum_{j=1}^{\infty}(2^j - 1 + \lambda)^{-k}$$

$$= \lambda^{-k}$$

and, hence,

$$\sum_{j=1}^{\infty}(2^j - 1 + \lambda)^{-k} = \int_0^{\infty} (2^x - 1 + \lambda)^{-k}\, dx + O(\lambda^{-k}).$$

For $k = 2$, this yields

(5.2)
$$\begin{aligned} S_\lambda^2 &= (\log 2)^{-1} \int_0^{\infty} (y + \lambda)^{-2}(y+1)^{-1}\, dy + O(\lambda^{-2}) \\ &= (\log \lambda)/(\lambda^2 \log 2) + O(\lambda^{-2}), \end{aligned}$$

as $\lambda \to \infty$. If we apply the same approach to $\lambda m_\lambda$ and $\log \varphi_J(\lambda)$, however, then the error caused by approximating these sums by integrals is of the order $O(1)$ and $O(\log \lambda)$, respectively, which yields a multiplicative error factor $(1 + O(a\lambda^b))$ in (5.1) for some positive $a$ and $b$. Of course this is not good enough so we shall have to expand the series representing $\lambda m_\lambda$ and $\log \varphi_J(\lambda)$ directly with remainder $o(1)$ in both cases.

Let $k$ be a natural number and $\theta \in [0, 1)$ be such that

$$\lambda = 2^{k+\theta},$$

and thus

$$k = (\log \lambda)/\log 2 - \theta = \lfloor (\log \lambda)/\log 2 \rfloor,$$

$$\theta = (\log \lambda)/\log 2 - k = \operatorname{frac}((\log \lambda)/\log 2).$$

Here $\lfloor x \rfloor$ and $\operatorname{frac}(x)$ are the integer and the fractional part of $x$, respectively. In order for $\lambda \to \infty$, it is necessary and sufficient that $k \to \infty$, while $\theta$ may vary arbitrarily in $[0, 1)$ with $k$. Using (5.2) we find

$$\begin{aligned} \lambda m_\lambda &= \sum_{j=1}^{\infty} \frac{\lambda}{2^j - 1 + \lambda} \\ &= \sum_{j=1}^{\infty} \frac{2^{k+\theta}}{2^j - 1 + 2^{k+\theta}} \\ &= \sum_{j=1}^{\infty} \frac{2^{k+\theta}}{2^j + 2^{k+\theta}} + O(\lambda^{-1} \log \lambda) \\ &= \sum_{j=1}^{k} \frac{1}{2^{j-k-\theta} + 1} + \sum_{j=k+1}^{\infty} \frac{1}{2^{j-k-\theta} + 1} + O(\lambda^{-1} \log \lambda) \end{aligned}$$

$$= \sum_{j=0}^{k-1} \frac{1}{2^{-j-\theta}+1} + \sum_{j=1}^{\infty} \frac{1}{2^{j-\theta}+1} + O(\lambda^{-1}\log\lambda)$$

$$= \sum_{j=1}^{k} \frac{2^j}{2^j+2^{1-\theta}} + \sum_{j=1}^{\infty} \frac{1}{2^{j-\theta}+1} + O(\lambda^{-1}\log\lambda)$$

$$= k - \sum_{j=1}^{k} \frac{2^{1-\theta}}{2^j+2^{1-\theta}} + \sum_{j=1}^{\infty} \frac{2^\theta}{2^j+2^\theta} + O(\lambda^{-1}\log\lambda)$$

$$= \frac{\log\lambda}{\log 2} - \sum_{j=1}^{\infty} \frac{2^{1-\theta}}{2^j+2^{1-\theta}} + \sum_{j=1}^{\infty} \frac{2^\theta}{2^j+2^\theta} - \theta + O(\lambda^{-1}\log\lambda).$$

Hence,

(5.3) $$\lambda m_\lambda = (\log\lambda)/\log 2 + A(\theta) + O(\lambda^{-1}\log\lambda),$$

with

$$A(\theta) = -\sum_{j=1}^{\infty} \frac{2^{1-\theta}}{2^j+2^{1-\theta}} + \sum_{j=1}^{\infty} \frac{2^\theta}{2^j+2^\theta} - \theta.$$

Notice that the term $A(\theta)$ of order 1 is not constant but depends on $\theta \in [0,1)$. The expansion (5.3) is obviously uniform in $\theta$.

Similarly, by (5.3),

$$\log \varphi_J(\lambda) = \sum_{j=1}^{\infty} \log\{(2^j-1)/(2^j-1+\lambda)\}$$

$$= \sum_{j=1}^{\infty} \log\{2^j/(2^j+\lambda)\} + \sum_{j=1}^{\infty} \log(1-2^{-j})$$

$$+ \sum_{j=1}^{\infty} \log\{1+1/(2^j-1+\lambda)\}$$

$$= \sum_{j=1}^{\infty} \log\{2^j/(2^j+2^{k+\theta})\} + \sum_{j=1}^{\infty} \log(1-2^{-j})$$

$$+ O(\lambda^{-1}\log\lambda).$$

Furthermore,

$$\sum_{j=1}^{\infty} \log\{2^j/(2^j+2^{k+\theta})\}$$

$$= \sum_{j=1}^{k} \log\{2^{1-\theta}/(2^j+2^{1-\theta})\} + \sum_{j=1}^{\infty} \log\{2^j/(2^j+2^\theta)\}$$

$$= k(1 - \theta) \log 2 - (1/2)k(k+1) \log 2$$
$$- \sum_{j=1}^{k} \log(1 + 2^{1-\theta-j}) - \sum_{j=1}^{\infty} \log(1 + 2^{\theta-j}).$$

Substituting $k = (\log \lambda)/\log 2 - \theta$ and using

$$\sum_{j=k+1}^{\infty} \log(1 + 2^{1-\theta-j}) \leq \sum_{j=k+1}^{\infty} 2^{1-\theta-j} = O(\lambda^{-1}),$$

we finally find

(5.4) $$\log \varphi_J(\lambda) = -\frac{(\log \lambda)^2}{2 \log 2} + \frac{\log \lambda}{2} + B(\theta) + O(\lambda^{-1} \log \lambda),$$

where

(5.5) $$B(\theta) = \sum_{j=1}^{\infty} \log\{(1 - 2^{-j})/[(1 + 2^{\theta-j})(1 + 2^{1-\theta-j})]\} - (1/2)\theta(1-\theta) \log 2.$$

Again the term $B(\theta)$ of order 1 depends on $\theta$ and the expansion is uniform in $\theta \in [0, 1)$.

Substituting (5.3), (5.4) and (5.2) in (5.1), we obtain, for $\lambda \to \infty$,

(5.6) $$P(m_\lambda) \sim \sqrt{\frac{\log 2}{2\pi \log \lambda}} \exp\left\{-\frac{1}{2 \log 2}(\log \lambda)^2 + \left(\frac{1}{2} + \frac{1}{\log 2}\right) \log \lambda + A(\theta) + B(\theta)\right\}.$$

It remains to find approximations of $\log \lambda$ and $(\log \lambda)^2$ as functions of

$$t = m_\lambda = (\log \lambda)/(\lambda \log 2) + A(\theta)/\lambda + O(\lambda^{-2} \log \lambda).$$

We find

$$\log(1/t) = \log \lambda - \log \log \lambda + \log \log 2 - A(\theta)(\log 2)/\log \lambda$$
$$+ O((\log \lambda)^{-2}),$$
$$\log \log(1/t) = \log \log \lambda - (\log \log \lambda)/\log \lambda + (\log \log 2)/\log \lambda$$
$$+ O((\log \log \lambda)^2/(\log \lambda)^2)$$

and, hence,

$$\log \lambda = \log(1/t) + \log \log(1/t) - \log \log 2 + (\log \log(1/t))/\log(1/t)$$
$$- (\log \log 2)/\log(1/t) + A(\theta)(\log 2)/\log(1/t)$$
$$+ O((\log \log(1/t))^2/(\log(1/t))^2),$$

$$(\log \lambda)^2 = [\log(1/t) + \log\log(1/t) - \log\log 2]^2 + 2\log\log(1/t)$$
$$- 2\log\log 2 + 2A(\theta)\log 2$$
$$+ O((\log\log(1/t))^2/\log(1/t)).$$

Together with (5.5) and (5.6), this yields that, for $t \to 0$,

$$P(t) \sim C(\theta) \exp\{-(2\log 2)^{-1}[\log(1/t) + \log\log(1/t) - \log\log 2]^2\}$$
$$\times t^{-(1/2+1/\log 2)}$$

with

(5.7) $$C(\theta) = 2^{-\theta(1-\theta)/2} \frac{1}{\sqrt{2\pi}} \prod_{j=1}^{\infty} \frac{1 - 2^{-j}}{(1 + 2^{\theta-j})(1 + 2^{1-\theta-j})}.$$

The factor $C(\theta)$ depends on $\theta = \text{frac}((\log \lambda)/\log 2)$.

It remains to express $\theta$ in terms of $t$. Define

$$\psi(t) = (\log 2)^{-1}[\log(1/t) + \log\log(1/t) - \log\log 2].$$

We have $k+\theta = (\log \lambda)/\log 2 = \psi(t)+o(1)$, and as $C$ is positive and bounded, the derivative of $C$ is positive and bounded and $C(\theta) = C(1-\theta)$. This implies that $C(\text{frac}\{\psi(t)\}) = C(\theta)(1 + o(1))$. It follows that, as $t \to +0$,

(5.8) $$P(t) \sim C(\text{frac}\{\psi(t)\}) \exp\left\{-\frac{\log 2}{2}[\psi(t)]^2\right\} t^{-(1/2+1/\log 2)},$$

with $C$ defined in (5.7). This is an exact asymptotic expression for $P(t)$ as $t \to +0$.

The dependence on $\text{frac}(\psi(t))$ in (5.8) is a most unusual feature. In fact, preliminary numerical calculations make one wonder whether there is any dependence at all, since one finds that $C(\theta)$ equals a constant ($\approx 0.01013$) throughout the interval $0 \leq \theta < 1$ to any reasonable degree of accuracy. Thus, in order to properly understand the asymptotic expression (5.8), we have to analyze $C(\theta)$ in more detail. Proposition 5.1 states that $C(\theta)$ does indeed depend on $\theta$, but in a very peculiar way. In fact, for any real $\theta$,

$$C(\theta) = \left[\frac{\sqrt{\log 2}}{2^{1/8} 2\pi} \prod_{j=1}^{\infty}(1 - 2^{-j})^2\right] (\tilde{\vartheta}_3(\theta))^{-1} \approx 0.01013(\tilde{\vartheta}_3(\theta))^{-1},$$

where

(5.9) $$\tilde{\vartheta}_3(\theta) = 1 + 2\sum_{k=1}^{\infty} \exp\{-2k^2\pi^2/\log 2\} \cos\{2k\pi(1/2 - \theta)\}$$
$$= \vartheta_3(\pi(1/2 - \theta), \exp\{-2\pi^2/\log 2\}).$$

Here $\vartheta_3$ is a theta function

$$\vartheta_3(z,q) = 1 + 2\sum_{k=1}^{\infty} q^{k^2} \cos(2kz).$$

Note that for all $\theta$,

$$|\tilde{\vartheta}_3(\theta) - 1| < 10^{-12}$$

is a quantity which is difficult to reveal numerically!

PROPOSITION 5.1. *For any real $\theta$,*

(5.10)
$$\prod_{j=1}^{\infty}(1+2^{\theta-j})(1+2^{1-\theta-j})$$
$$= 2^{-\theta(1-\theta)/2}\,\tilde{\vartheta}_3(\theta)\frac{2^{1/8}\sqrt{2\pi}}{\sqrt{\log 2}}\prod_{j=1}^{\infty}(1-2^{-j})^{-1},$$

*where $\tilde{\vartheta}_3(\theta)$ is given by* (5.9).

PROOF. We first apply Jacobi's triple product identity [see, e.g., Askey (1980) and Gasper and Rahman (1990)]. For any $q \in (0,1)$,

$$(5.11) \quad \prod_{j=0}^{\infty}(1-xq^j)(1-x^{-1}q^{j+1})(1-q^{j+1}) = \sum_{n=-\infty}^{\infty}(-1)^n q^{\binom{n}{2}}x^n.$$

Take $x = -2^{-\theta}$, $q = 1/2$. Then (5.11) becomes

$$(5.12) \quad \prod_{j=1}^{\infty}(1+2^{\theta-j})(1+2^{1-\theta-j})(1-2^{-j}) = \sum_{n=-\infty}^{\infty} 2^{-n(n-1)/2}2^{-\theta n}.$$

The right-hand side of (5.12) is of the form

$$c(\theta)\sum_{n=-\infty}^{\infty} g(n),$$

where

$$c(\theta) = \frac{2^{1/8}\sqrt{2\pi}}{\sqrt{\log 2}} 2^{-\theta(1-\theta)/2},$$

and

$$g(x) = \frac{\sqrt{\log 2}}{\sqrt{2\pi}} \exp\{-(1/2)(\log 2)(x - 1/2 + \theta)^2\}$$

is a normal density with mean $\mu = 1/2 - \theta$ and standard deviation $\sigma = 1/\sqrt{\log 2}$. The characteristic function of $g$ is given by

$$\gamma(t) = \exp\{-t^2/(2\log 2) + it(1/2 - \theta)\},$$

where $i$ is the imaginary unit. For each fixed $\lambda$ and for each real $\xi$, the Poisson summation formula [see Feller (1970)] gives

$$(5.13) \qquad \sum_{k=-\infty}^{+\infty} \gamma(\xi + 2k\lambda) = \frac{\pi}{\lambda} \sum_{n=-\infty}^{+\infty} g(n\pi/\lambda) \exp\{in(\pi/\lambda)\xi\}.$$

Put $\lambda = \pi$, $\xi = 0$. Then the right-hand side of (5.13) becomes $\sum_{n=-\infty}^{\infty} g(n)$ and on the left-hand side we have

$$\sum_{k=-\infty}^{\infty} \gamma(2k\pi) = \gamma(0) + \sum_{\substack{k=-\infty \\ k \neq 0}}^{\infty} \exp\{-2k^2\pi^2/\log 2\} \exp\{i(1/2 - \theta)2k\pi\}$$

$$= 1 + 2\sum_{k=1}^{\infty} \exp\{-2k^2\pi^2/\log 2\} \cos\{2k\pi(1/2 - \theta)\}$$

$$= \tilde{\vartheta}_3(\theta).$$

Hence, (5.13) reduces to

$$\sum_{n=-\infty}^{+\infty} g(n) = \tilde{\vartheta}_3(\theta),$$

implying that the right-hand side of (5.12) equals $c(\theta)\tilde{\vartheta}_3(\theta)$. This immediately yields (5.10). The proposition is proved. $\square$

We summarize our findings in the following theorem.

THEOREM 5.2. *Let $X_1, X_2, \ldots$ be independent exponential random variables with mean 1, and let*

$$J = \sum_{j=1}^{\infty} (2^j - 1)^{-1} X_j.$$

*Then*

$$\mathbb{P}(J \leq t) \sim \frac{\sqrt{\log 2} \prod_{j=1}^{\infty}(1 - 2^{-j})^2}{2^{1/8} 2\pi \tilde{\vartheta}_3(\mathrm{frac}\{\psi(t)\})}$$

$$\times \exp\left\{-\frac{\log 2}{2}[\psi(t)]^2\right\} t^{-(1/2 + 1/\log 2)} \qquad as\ t \to +0,$$

*where*

$$\psi(t) = (\log 2)^{-1}[\log(1/t) + \log\log(1/t) - \log\log 2]$$

*and $\tilde{\vartheta}_3$ is defined in (5.9).*

**6. Related results.** In a similar fashion we can also analyze more general linear combinations of i.i.d. exponential random variables than $J$. For any $q \in (0,1)$, define

$$J^{(q)} = (q^{-1} - 1) \sum_{j=1}^{\infty} (q^{-j} - 1)^{-1} X_j.$$

Clearly, $J \equiv J^{(1/2)}$. One can show that

$$(m+1)(q^{-1} - 1) \sum_{j=1}^{m+1} \frac{1}{q^{-j} - 1} D_{j,m} \xrightarrow{d} J^{(q)} \qquad \text{as } m \to \infty,$$

where the expression on the left occurs in the right-hand side of (1.5) for $n = m$. The random variable $J^{(q)}$ can be written in the following way. Let $N(t)$ be a standard Poisson process. Then

$$(q^{-1} - 1)^{-1} J^{(q)} = \int_0^{\infty} q^{N(t)+1} (1 - q^{N(t)+1})^{-1} \, dt$$

$$= \sum_{j=1}^{\infty} q^j \int_0^{\infty} q^{jN(t)} \, dt$$

$$= \sum_{j=1}^{\infty} q^j I^{(q^j)},$$

where

$$I^{(q)} = \int_0^{\infty} q^{N(t)} \, dt = \sum_{j=1}^{\infty} q^{j-1} X_j$$

is an exponential functional associated with a Poisson process. The functional $I^{(q)}$ has been intensively studied in recent literature. Its density was obtained independently by Dumas, Guillemin and Robert (2002), Bertoin, Biane and Yor (2002) and Litvak and Adan (2001), for $q = 1/2$. Carmona, Petit and Yor (1997) derived a density of $\int_0^{\infty} h(N(t)) \, dt$ for a large class of functions $h : \mathbb{N} \to \mathbb{R}_+$, in particular, for $h(n) = q^n$. Bertoin, Biane and Yor (2002) found the fractional moments of $I^{(q)}$. If $i^{(q)}(t)$ is a density of $I^{(q)}$, then $i^{(q)}(t)$ and all its derivatives equal 0 at the point $t = 0$. This implies, by the way, that all moments of $1/I^{(q)}$ are finite. However, for $q = 1/e$, it was proved by Bertoin and Yor (2002a) that $1/I^{(1/e)}$ is not determined by its moments.

The functional $I^{(q)}$ appears in a number of applications. Let $T_n^{NI}$ be the travel time needed to collect $n$ items independently and uniformly distributed on a circle of length 1 operating under the nearest item heuristic (the picker always travels to the nearest item to be retrieved). Then it was shown by Litvak and Adan (2001) that $(n+1)(1 - T_n^{NI})$ converges in distribution to $I^{(1/2)}$. Dumas, Guillemin and Robert (2002) showed that the distribution of $I^{(q)}$ plays a key role in the analysis of limiting behavior of a Transmission Control Protocol connection. These results were extended by Guillemin, Robert and Zwart (2002), who found the distribution and the fractional moments of the exponential functional

$$(6.1) \qquad I(\xi) = \int_0^\infty e^{-\xi(t)}\, dt,$$

where $(\xi(t), t \geq 0)$ is a compound Poisson process. An exponential functional (6.1) associated with a Levy process $\xi(t)$ appears in mathematical finance and many other fields. It has been studied recently by Bertoin and Yor (2001, 2002a, b), Bertoin, Biane and Yor (2002), Carmona, Petit and Yor (1997) and Yor (2001).

Along the same lines as in Section 5, one can prove theorems similar to Theorem 5.2 for $I^{(q)}$ and $J^{(q)}$. In fact, it is straightforward to repeat the calculations for

$$qI^{(q)} = \sum_{j=1}^\infty q^j X_j \quad \text{and} \quad \frac{q}{1-q} J^{(q)} = \sum_{j=1}^\infty \frac{1}{q^{-j}-1} X_j.$$

We obtain, for $q \in (0,1)$ as $t \to +0$,

$$\mathbb{P}(qI^{(q)} \leq t) \sim \frac{1}{2\pi} q^{1/8} \sqrt{\log(1/q)} \left[\prod_{j=1}^\infty (1-q^j)\right] t^{-(1/2+1/\log(1/q))}$$

$$\times \exp\left\{-\frac{\log(1/q)}{2} [\psi^{(q)}(t)]^2\right\} [\tilde{\vartheta}_3^{(q)}(\text{frac}\{\psi^{(q)}(t)\})]^{-1},$$

$$\mathbb{P}\left(\frac{q}{1-q} J^{(q)} \leq t\right) \sim \mathbb{P}(qI^{(q)} \leq t) \prod_{j=1}^\infty (1-q^j),$$

where

$$\psi^{(q)}(t) = (\log(1/q))^{-1}[\log(1/t) + \log\log(1/t) - \log(\log(1/q))],$$

$$\tilde{\vartheta}_3^{(q)}(\theta) = 1 + 2\sum_{k=1}^\infty \exp\{-2k^2\pi^2/\log(1/q)\} \cos\{2k\pi(1/2-\theta)\}.$$

This agrees with the result of Bertoin and Yor (2002a) that

$$\log i(t) \sim -\tfrac{1}{2}(\log(1/t))^2 \qquad \text{as } t \to +0,$$

where $i(t)$ is a density of

$$I = \int_0^\infty e^{-N(t)}\, dt = \sum_{j=1}^\infty e^{-j} X_j.$$

For the functional $I^{(1/2)}$, which describes the limiting behavior of the travel time under the nearest item heuristic, we find

$$\mathbb{P}(I^{(1/2)} \leq 2t) \sim \mathbb{P}(J \leq t) \prod_{j=1}^{\infty} (1 - 2^{-j})^{-1}, \qquad t \to +0.$$

**Acknowledgments.** We are grateful to Bert Zwart for drawing our attention to relevant literature and particularly to the paper of Davis and Resnick (1991), which inspired the results in Section 5. Discussions with Fred Steutel also proved very helpful.

FACULTY OF ELECTRICAL ENGINEERING  
MATHEMATICS AND COMPUTER SCIENCE  
UNIVERSITY OF TWENTE  
P.O. BOX 217  
7500 AE ENSCHEDE  
THE NETHERLANDS  
E-MAIL: n.litvak@math.utwente.nl

DEPARTMENT OF MATHEMATICS  
UNIVERSITY OF LEIDEN  
P.O. BOX 9512  
2300 RA LEIDEN  
THE NETHERLANDS  
E-MAIL: vanzwet@math.leidenuniv.nl